\documentclass[a4paper,11pt]{article}
\usepackage[utf8]{inputenc}
\usepackage{color}

\usepackage[pagewise]{lineno}

\usepackage[all]{xy}
\usepackage{amsmath,amssymb,amsfonts,amsthm,graphicx}
\usepackage[margin=1.15in]{geometry}
\usepackage{multirow}
\theoremstyle{plain}

\usepackage{tikz}
\usetikzlibrary{matrix,arrows,positioning,shapes.geometric}

\tikzset{node distance=2cm, auto}

\newtheorem{thm}{Theorem}[section]

\newtheorem{cor}[thm]{Corollary}
\newtheorem{prop}[thm]{Proposition}
\newtheorem{df}[thm]{Definition}
\newtheorem{lm}[thm]{Lemma}

\newtheorem*{theorem*}{Theorem}
\newtheorem*{aim*}{Aim}
\newtheorem*{initialaim*}{Initial Aim}
\newtheorem*{conj*}{Conjecture}
\newtheorem*{cor*}{Corollary}
\newtheorem*{prop*}{Proposition}
\newtheorem*{df*}{Definition}
\newtheorem*{lm*}{Lemma}
\newtheorem*{example*}{Example}
\newtheorem*{notation*}{Notation}
\newtheorem*{prob*}{Problem}

\numberwithin{thm}{section}

\title{Central intersections of element centralisers}
\author{Julian Brough}

\begin{document}
\date{}
\maketitle

\vspace{-3mm}

\begin{center}
\small
\textit{FB Mathematik, TU Kaiserslautern, Postfach 3049}

\textit{67653 Kaiserslautern, Germany}

\text{E-mail: brough@mathematik.uni-kl.de}
\end{center}

\paragraph{}

\textit{Seconday: 20D99}

\paragraph{}

\normalsize
\begin{abstract}
In 1970 R. Schmidt gave a structural classification for CA-groups.
In this paper we consider a condition upon the intersection of element centralisers which turns out to be equivalent to the definition of a CA-group.
We then weaken which centralisers we chose to intersect and structurally classify this new family of groups.
Furthermore we apply a similar weakening to the class of F-groups introduced by It{\^o} in 1953 and classified by Rebmann in 1971.
 
\end{abstract}

\section{Introduction}

A finite group is called a CA-group if the centraliser of every non-central element is abelian.
If $G$ is a CA-group and $x,y\in G\setminus Z(G)$, then $C_G(x)$ can never be properly contained in $C_G(y)$.
Furthermore we shall see in Lemma~\ref{TIC=CA} that $G$ being a CA-group is equivalent to saying for all $x$ and $y$ non-central elements in $G$ either $C_G(x)=C_G(y)$ or $C_G(x)\cap C_G(y)=Z(G)$.
In addition to the class of CA-groups, It{\^o} \cite{ItoTypeI} introduced the notion of an F-group.
This is a group in which every non-central element centraliser contains no other non-central element centraliser.
That is, $G$ is an F-group if for any $x\in G\setminus Z(G)$, then $C_G(x)\leq C_G(y)$ implies $y\in Z(G)$.
We shall also see in Lemma~\ref{CentralTIC=F} that $G$ being an F-group is equivalent to $Z(C_G(x))\cap Z(C_G(y))=Z(G)$ for all $x,y\in G\setminus Z(G)$ such that $C_G(x)\ne C_G(y)$. 

The aim of this paper is to consider these intersection conditions for a specific subset of centralisers, in particular, the set of minimal centralisers in a group (those which do not properly contain any other element centraliser).
Thus we define a group to be a ${\rm CA}_{min}$-group if for two non-central elements $x$ and $y$ with minimal element centralisers in $G$ either $C_G(x)=C_G(y)$ or $C_G(x)\cap C_G(y)=Z(G)$.
Similarly we call $G$ an ${\rm F}_{min}$-group if for two non-central elements $x$ and $y$ with minimal element centralisers in $G$ either $C_G(x)=C_G(y)$ or $Z(C_G(x))\cap Z(C_G(y))=Z(G)$.
Note that the analogous condition by considering maximal centralisers for CA-groups was studied by Schmidt \cite{SchmidtCaGps} (although Schmidt defined such groups using subgroup centralisers).

The aim of this paper is to prove a structural classification of ${\rm CA}_{min}$-groups and ${\rm F}_{min}$-groups.

\begin{thm}\label{MainThm}
$G$ is a ${\rm CA}_{min}$-group {\rm (respectively ${\rm F}_{min}$)} if and only if $G$ has one of the following forms:
\begin{enumerate}
\item $G/Z(G)$ is a Frobenius group with kernel $L/Z(G)$ and complement $K/Z(G)$ such that both $K$ and $L$ are abelian.
\item $G/Z(G)$ is a Frobenius group with kernel $L/Z(G)$ and complement $K/Z(G)$ such that $K$ is abelian, $L$ is a ${\rm CA}_{min}$-group {\rm (respectively ${\rm F}_{min}$)}, $Z(G)=Z(L)$ and $L/Z(L)$ is a $p$-group.
\item $G/Z(G)\cong {\rm Sym}(4)$ and if $V/Z(G)\cong V_4$, then $V$ is non-abelian. 
\item $G$ has an abelian normal subgroup of index $p$, $G$ is not abelian.
\item $G\cong A\times P$, where $A$ is abelian and $P$ is a non-abelian $p$-group for some prime $p$; therefore $P$ is a ${\rm CA}_{min}$-group {\rm (respectively ${\rm F}_{min}$)}.
\item $G/Z(G)\cong PGL_2(p^n)$ or $PSL_2(p^n)$ with $p^n > 3$.
\end{enumerate}
\end{thm}

Note a similarity to Rebmann's structural classification of F-groups \cite{FGroups}.
In fact the groups of type (1), (3) and (4) are CA-groups \cite{SchmidtCaGps}.
Moreover for the families (2) and (5) replacing ${\rm CA}_{min}$-group by F-group yields the corresponding family for F-groups.
While the non-solvable case (6) contains all the non-solvable cases of F-groups.
In particular we obtain the following corollary as given in Rebmann \cite{FGroups}.

\begin{cor}\label{NonSolFIsCA}
Let $G$ be a non-solvable F-group.
Then $G$ is a CA-group.
\end{cor}

By using the above theorem is it easy to see that the class of CA-groups is  strictly smaller than the class of ${\rm CA}_{min}$-groups, as $PSL_2(q)$ and $PGL_2(q)$ have a non-abelian centraliser.
However in Rebmann's paper no example of an F-group which is not a CA-group was provided.
We finish this paper by providing a family of $p$-groups which are F-groups but not CA-groups.

\begin{prop}\label{FNotCA}
Let $G$ be an extraspecial group of order $p^{2n+1}$ with $n>1$.
Then $G$ is an F-group which is not a CA-group.
\end{prop}

Note that this family will also be a family of ${\rm F}_{min}$-groups which are not ${\rm CA}_{min}$-groups.
Finally, observe that if there exists a solvable ${\rm CA}_{min}$-group which is not an CA-group, then such a $p$-group exists for some prime $p$. 
However, running over the GAP libraries we have been unable to find a $p$-group which is ${\rm CA}_{min}$-group and not a CA-group.
Note that \cite{RockeAbCent} studied such groups however we were unable to use the results and methods in this paper to produce such an example.

\section{Preliminaries}
\subsection{Conditions on element centraliser intersections}

Let $G$ be a non-abelian CA-group with $x$ and $y$ non-central elements in $G$.
Consider $C_G(x)\cap C_G(y)$.
If $z\in C_G(x)\cap C_G(y)$, then $\langle C_G(x),C_G(y)\rangle\leq C_G(z)$.
As these centralisers are abelian either $C_G(x)= C_G(y)$ or $z\in Z(G)$.
Moreover we shall show in the next lemma that this condition is equivalent to a CA-group. 

\begin{lm}\label{TIC=CA}
Let $G$ be a finite non-abelian group.
Then $G$ is a CA-group if and only if for any pair of non-central elements $x$ and $y$ such that $C_G(x)\ne C_G(y)$ then $C_G(x)\cap C_G(y)=Z(G)$.
\begin{proof}
We commented above that any CA-group satisfies this condition. 
Hence it remains to show to converse.
Let $z$ be a non-central element in $G$ and $x,y\in C_G(z)\setminus Z(G)$.
Then $\langle z\rangle \leq C_G(y)\cap C_G(x)$ and so $C_G(x)\cap C_G(y)\ne Z(G)$.
Therefore $C_G(x)= C_G(y)$, and so $x$ and $y$ commute.
In particular, any two elements in $C_G(z)$ commute.
Thus $C_G(z)$ is abelian and $G$ is a CA-group.
\end{proof}
\end{lm}

Furthermore, with this observation we have the following corollary.
\begin{cor}
Let $G$ be a CA-group such that $Z_2(G)>Z(G)$.
Then $G$ is meta-abelian.
\begin{proof}
We want to show that $G'$ is abelian.
However, by \cite[Theorem III.2.11]{Huppert}, $G'\leq C_G(Z_2(G))$.
Thus it is enough to show $C_G(Z_2(G))$ is abelian.
Let $x,y\in C_G(Z_2(G))\setminus Z(G)$.
Then $Z(G)<Z_2(G)\leq C_G(x)\cap C_G(y)$.
Thus $x$ and $y$ commute and so $C_G(Z_2(G))$ is abelian.
 
\end{proof}
\end{cor}

We now make clear the definition of a minimal element centraliser for use in the definitions of ${\rm CA}_{min}$-groups and ${\rm F}_{min}$-groups.

\begin{df}
An element centraliser $C_G(x)$ for $x$ a non-central element is called a minimal centraliser if $C_G(y)\leq C_G(x)$ implies $C_G(y)=C_G(x)$.
\end{df}

Thus to relax the notion of a CA-group, we want to consider the intersection property for minimal element centralisers.
Note that any non-abelian group must have at least two minimal non-central element centralisers.
Otherwise if $C=C_G(x)$ is the unique minimal centraliser in $G$, then for all $y\in G$, we have that $C\leq C_G(y)$.
Therefore $x\in \cap_{y\in G} C_G(y)=Z(G)$.
Thus $C=G$ and $x\in Z(G)$.

Note that we could also consider maximal centralisers ($C_G(x)$ called a maximal centraliser if $C_G(x)<C_G(y)$ implies $y\in Z(G)$).
In fact Schmidt considered the set of groups in which any two distinct maximal non-central element centralisers have intersection equal to the center of the group \cite{SchmidtCaGps}. (These were referred to as $\mathfrak{D}$-groups)
However he only classified the soluble $\mathfrak{D}$-groups, although he did discuss in depth the non-solvable case too.

\begin{df}
Let ${\rm CA}_{min}$ denote the set of finite groups $G$ such that $C_G(x)\cap C_G(y)=Z(G)$ for any two distinct minimal centralisers $C_G(x)$ and $C_G(y)$. 
\end{df}

By definition, a group is an F-group if and only if for any non-central element $x$ in $G$ we have $C_G(x)$ is both a maximal and minimal centraliser in $G$. 
Therefore, for an F-group, the definitions of $\mathfrak{D}$ and ${\rm CA}_{min}$ are equivalent.
Furthermore, the following corollary follows from Lemma~\ref{TIC=CA}.

\begin{cor}\label{CA=D+F}
Let $G$ be a finite group.
Then $G$ is a CA-group if and only if $G$ is an F-group and a ${\rm CA}_{min}$-group.
\end{cor}
 
As with the notion of CA-groups we shall weaken the notion of an F-group.
However, first we need an analogous lemma for Lemma~\ref{TIC=CA}.

\begin{lm}\label{CentralTIC=F}
Let $G$ be a finite non-abelian group.
Then $G$ is an F-group if and only if for any pair of non-central elements $x$ and $y$ such that $C_G(x)\ne C_G(y)$ then $Z(C_G(x))\cap Z(C_G(y))=Z(G)$.
\begin{proof}

Assume $G$ is an F-group and let $C_G(x)\ne C_G(y)$ for $x$ and $y$ non-central elements.
If $z\in Z(C_G(x))\cap Z(C_G(y))$, then $\langle C_G(x),C_G(y)\rangle \leq C_G(z)$. 
Hence as in an F-group every centraliser is both maximal and minimal, it follows that $C_G(z)=G$.
Or in other words $z\in Z(G)$.
Therefore $Z(C_G(x))\cap Z(C_G(y))=Z(G)$.

For the converse direction assume that $C_G(x)<C_G(y)$ for both $x$ and $y$ non-central elements in $G$.
Then $Z(C_G(y))\leq Z(C_G(x))$ and therefore $Z(C_G(y))=Z(G)$, which implies $y\in Z(G)$.
\end{proof}
\end{lm}

Thus we now make the following definition.

\begin{df}
Let ${\rm F}_{min}$ denote the set of finite groups $G$ such that $Z(C_G(x))\cap Z(C_G(y))=Z(G)$ for any two distinct minimal centralisers $C_G(x)$ and $C_G(y)$. 
\end{df}

Therefore we have the following inclusions:
(Theorem~\ref{MainThm} and Proposition~\ref{FNotCA} show that these are strict inclusions) 

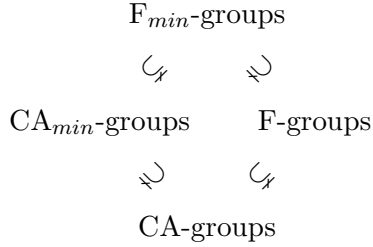
\begin{center}
\begin{tikzpicture}
\node (1) {${\rm F}_{min}$-groups};
\node[below left of=1,node distance=10mm , rotate=45] (2) {$\subsetneq$};
\node[below right of=1,node distance=10mm , rotate=315] (3) {$\supsetneq$};
\node[below left of=1] (4) {${\rm CA}_{min}$-groups};
\node[below right of=1] (5) {F-groups};
\node[below left of=5] (6) {CA-groups};
\node[above right of=6,node distance=10mm , rotate=45] (7) {$\subsetneq$};
\node[above left of=6,node distance=10mm , rotate=315] (8) {$\supsetneq$};

\end{tikzpicture}
\end{center}

We finally observe that the intersection of ${\rm CA}_{min}$-groups with F-groups equals the set of CA-groups. 

\subsection{Exhibiting a partition}

In the works of Rebmann and Schmidt \cite{FGroups}, \cite{SchmidtCaGps} providing an abelian normal partition of the central quotient $G/Z(G)$ yielded a powerful tool to structurally classify families of groups; in particular they could apply the following classifications by Baer and Suzuki.
Neither theorem appears as one statement but as several across the papers, therefore we combine the results into one statement.

\begin{thm}\cite{BaerPart1}\cite{BaerPart2}\label{BaerSolPart}
Let $G$ be a solvable group with a normal non-trivial partition $\beta$, then $G$ is one of the following:
\begin{enumerate}
\item A component of $\beta$ is self normalising in $G$ and $G$ is a Frobenius group.
\item $G\cong {\rm Sym}(4)$ and $\beta$ is the set of maximal cyclic subgroups of $G$.
\item $G$ has a nilpotent normal subgroup $N$ which lies in $\beta$ with $|G:N|=p$ and every element in $G\setminus N$ has order $p$.
\item $G$ is a $p$-group, for $p$ a prime.
\end{enumerate}
\end{thm}

\begin{thm}\cite{SuzPart}\label{SuzNonSolPart}
Let $G$ be a non-solvable group with a normal non-trivial partition $\beta$.
Then $G\cong PGL_2(p^n)$, $PSL_2(p^n)$ for $p$ prime and $p^n>3$, $Sz(2^n)$ for $n\geq 3$ or a component of $\beta$ is self normalising and $G$ is a Frobenius group.
\end{thm}

We aim to show that for $G$ a ${\rm CA}_{min}$-group or an ${\rm F}_{min}$-group, as for F-groups, the central quotient $G/Z(G)$ exhibits a normal abelian partition.
For the case of ${\rm CA}_{min}$-groups we require the following preliminary result.

\begin{lm}\label{CAminCenAb}
Let $G$ be a ${\rm CA}_{min}$-group, then each minimal centraliser is abelian.
\begin{proof}
Let $C$ be a minimal centraliser in $G$ and $x\in C$.
If $x\not\in Z(C)$, then there exists a minimal centraliser $D \leq C_G(x)$.
It is clear that $Z(C_G(x))\leq Z(D)$.
Hence $x\in C\cap D$ which equals $Z(G)$ or $C=D$.
Thus assume that $C=D$.
However as $x\in Z(D)$, it means that $x\in Z(C)$.
\end{proof}
\end{lm}

\begin{lm}
Let $G$ be a ${\rm CA}_{min}$-group.
Then 
\[
\beta  = \{C/Z(G) \mid C \text{ a minimal centraliser in $G$}\}
\]
forms a non-trivial normal partition of $G/Z(G)$ consisting of abelian subgroups.
\begin{proof}
It is clear that the set $\beta$ is closed under conjugation and by Lemma~\ref{CAminCenAb} every subgroup in $\beta$ is abelian.
Thus to show $\beta$ is a partition we need to show that every element in $G/Z(G)$ lies in a unique subgroup in $\beta$.

Take $C/Z(G)$ and $D/Z(G)$ distinct in $\beta$.
Then $C/Z(G)\cap D/Z(G)=(C\cap D)/Z(G)=1$.
Thus it is enough to show that any $xZ(G)$ lies in some $C/Z(G)$.

Consider $C_G(x)$ which contains some minimal centraliser $C$. 
Then as in Lemma~\ref{CAminCenAb}, $Z(C_G(x))\leq Z(C)$, hence $x\in C$.
In particular $xZ(G)\in C/Z(G)$.
\end{proof}
\end{lm}

\begin{lm}
Let $G$ be an ${\rm F}_{min}$-group.
Then 
\[
\beta  = \{Z(C)/Z(G) \mid C \text{ a minimal centraliser in $G$}\}
\]
forms a non-trivial normal partition of $G/Z(G)$ consisting of abelian subgroups.
\begin{proof}
It is clear that the set $\beta$ is closed under conjugation and every subgroup in $\beta$ is abelian.
Thus to show $\beta$ is a partition we need to show that every element in $G/Z(G)$ lies in a unique subgroup in $\beta$.

Take $Z(C)/Z(G)$ and $Z(D)/Z(G)$ distinct in $\beta$.
Then $Z(C)/Z(G)\cap Z(D)/Z(G)=(Z(C)\cap Z(D))/Z(G)=1$.
Thus it is enough to show that any $xZ(G)$ lies in some $Z(C)/Z(G)$.

Consider $C_G(x)$ which contains some minimal centraliser $C$. 
Then as in Lemma~\ref{CAminCenAb}, $Z(C_G(x))\leq Z(C)$, hence $x\in Z(C)$.
In particular $xZ(G)\in Z(C)/Z(G)$.
\end{proof}
\end{lm}

\section{Proof of main theorem}

We now aim to classify ${\rm CA}_{min}$-groups and ${\rm F}_{min}$-groups.
In fact a similar argument as for F-groups occurs when we replace F-group by ${\rm CA}_{min}$-group or ${\rm F}_{min}$-group.

\subsection{Classifying ${\rm CA}_{min}$-groups}

Due to the classification of partitions by Baer and Suzuki, first we shall consider the solvable ${\rm CA}_{min}$-groups.

\begin{thm}\label{SolCAminGp}
Let $G$ be a solvable ${\rm CA}_{min}$ group.
Then $G$ is one of the following:
\begin{enumerate}
\item $G/Z(G)$ is a Frobenius group with kernel $L/Z(G)$ and complement $K/Z(G)$ such that both $K$ and $L$ are abelian.
\item $G/Z(G)$ is a Frobenius group with kernel $L/Z(G)$ and complement $K/Z(G)$ such that $K$ is abelian, $L$ is a ${\rm CA}_{min}$-group, $Z(L)=Z(G)$ and $L/Z(L)$ is a $p$-group.
\item $G/Z(G)\cong {\rm Sym}(4)$ and if $V/Z(G)\cong V_4$, then $V$ is non-abelian. 
\item $G$ has an abelian normal subgroup of index $p$, $G$ is not abelian.
\item $G\cong A\times P$, where $A$ is abelian and $P$ is a non-abelian $p$-group for some prime $p$; therefore $P$ is a ${\rm CA}_{min}$-group.
\end{enumerate}

\begin{proof}
As $G/Z(G)$ admits a non-trivial normal partition, we apply the classification of Baer (Theorem~\ref{BaerSolPart}) to determine $G/Z(G)$.

\underline{\bf Case (1)}\newline
Let $L/Z(G)$ denote the Frobenius kernel of $G/Z(G)$.
Furthermore, let $K/Z(G)$ denote an element in the partition of $G/Z(G)$ which is self-normalising.
Then $K=C_G(x)$ for some minimal centraliser $C_G(x)$ in $G$.
We want to show that $K/Z(G)$ is a Frobenius complement, thus we need that $K\cap K^g=Z(G)$ for all $g\in G\setminus K$.

As $K$ is a minimal centraliser in $G$, it means  $K^g$ is also a minimal centraliser in $G$ and therefore $K\cap K^g=Z(G)$ or $K$.
However, if $K^g=K$, then $gZ(G)\in N_{G/Z(G)}(K/Z(G))=K/Z(G)$ as $K/Z(G)$ was chosen to be self-normalising.
Thus $K/Z(G)$ is a Frobenius complement in $G/Z(G)$.
Furthermore as $K$ is a minimal centraliser in $G$, then $K$ is abelian (Corollary~\ref{CAminCenAb}).

Let $x\in L\setminus Z(G)$.
As $G/Z(G)$ is Frobenius with kernel $L/Z(G)$, 
\[
C_G(x)/Z(G)\leq C_{G/Z(G)}(xZ(G))\leq L/Z(G),
\]
so $C_G(x)\leq L$; or in other words $C_G(x)=C_L(x)$.

If $L$ has a unique minimal centraliser, then $L$ is abelian by Corollary~\ref{CAminCenAb} and thus of type $(1)$.
Hence assume $L$ has two distinct minimal centralisers $C_L(x)<L$ and $C_L(y)<L$ for $x,y\in L\setminus Z(G)$.
Then $C_L(x)=C_G(x)$ and $C_L(y)=C_G(x)$.
If $C_G(x)$ is not a minimal centraliser in $G$, then there exists $C_G(z)<C_G(x)$.
As $C_G(z)<C_G(x)\leq L$, it follows that $z\in L\setminus Z(G)$.
Thus $C_L(z)< C_L(x)$ and so $C_L(x)$ is not minimal.
Therefore $C_G(x)$ and $C_G(y)$ are distinct minimal centralisers in $G$.
As $G$ is a ${\rm CA}_{min}$-group, we have that $Z(L)\leq C_L(x)\cap C_L(y)=C_G(x)\cap C_G(y)=Z(G)$.
However $Z(G)\leq L$ and therefore $Z(G)\leq Z(L)$ implying that $Z(L)=Z(G)$.
Furthermore, we have shown that $L$ is a ${\rm CA}_{min}$-group.
By Thompson, \cite[Theorem V.8.7]{Huppert}, $L/Z(G)$ is nilpotent (as it is a Frobenius kernel).
By \cite[Remark 2.4]{BaerPart1}, the only nilpotent groups with a partition are $p$-groups for some prime $p$.
This implies $G$ is of type $(2)$.

\underline{\bf Case (2)}\newline
In this case $G/Z(G)\cong {\rm Sym}(4)$. 
Let $V\leq G$ such that $V/Z(G)\cong V_4$, the Klein-four subgroup.
If $V$ is abelian, then for all $x\in V\setminus Z(G)$ we have $V\leq C_G(x)$.
As $xZ(G)$ has order $2$ it follows that $C_{G/Z(G)}(xZ(G))=D_8$ or $V_4$ (when $xZ(G)$ is a double or single transposition respectively), we also observe that $V/Z(G)\leq C_G(x)/Z(G)\leq  C_{G/Z(G)}(xZ(G))$.

If there exists an $x\in V$ such that $C_G(x)/Z(G)\cong V_4$, then $C_G(x)=V$ and so is abelian and thus a minimal centraliser.
However $V_4$ is not contained in the partition of $G/Z(G)$ which yields a contradiction. 
Thus for all $x\in V\setminus Z(G)$ we have that $C_G(x)/Z(G)=C_{G/Z(G)}(xZ(G))\cong D_8$ and $xZ(G)\leq Z(D_8)$.
In particular, it follows that $V_4$ must be the normal klein-four subgroup of ${\rm Sym}(4)$.

Inside $C_G(x)/Z(G)\cong D_8$ there exists a unique cyclic subgroup of order 4 and another copy of $V_4$ which is not normal in ${\rm Sym}(4)$.
Let $N,M\leq C_G(x)$ such that $N/Z(G)\cong C_4$ and $M/Z(G)$ equals the non-normal copy of $V_4$.
By Theorem~\ref{BaerSolPart}, $N/Z(G)$ lies in the partition and therefore $N$ is a minimal centraliser in $G$.
The subgroup $M$ contains $x$ and therefore $\langle x,Z(G)\rangle \leq Z(M)$. 
In particular, we see that $M/Z(M)$ must be cyclic and so $M$ is abelian.
However $x\in N\cap M$ and so $M$ cannot be a minimal centraliser in $G$.
Thus $C_G(y)>M$ for all $y\in M$.
As $M/Z(G)$ lies in a unique maximal subgroup isomorphic to $D_8$, it follows that for each $y\in M\setminus Z(G)$, then $C_G(y)/Z(G)$ equals the unique maximal subgroup containing $M/Z(G)$.
Thus $M/Z(G)\leq Z(D_8)$ which is a contradiction.

\underline{\bf Case (3)}\newline
In this case $N/Z$ is a component of $\beta$, which implies $N$ is abelian.

\underline{\bf Case (4)}\newline
In this case $G/Z(G)$ is a $p$-group and therefore $G$ is nilpotent.
Therefore $G=A\times P$ for $A\leq Z(G)$ and $P$ a $p$-subgroup which is a ${\rm CA}_{min}$-group.
\end{proof}

\end{thm}

We next show that each case occurring in Theorem~\ref{SolCAminGp} yields a ${\rm CA}_{min}$-group. 

\begin{prop}\label{ListAreCAmin}
Any solvable group occurring in Theorem~\ref{SolCAminGp} is a ${\rm CA}_{min}$-group.
\begin{proof}
The solvable groups in Theorem~\ref{MainThm} are those of type $(1)-(5)$.
Any group of type $(5)$ is easily seen to be a ${\rm CA}_{min}$-group.
For the groups of type $(1),(3)$ and $(4)$, Schmidt \cite{SchmidtCaGps} has shown that they are CA-groups and hence are ${\rm CA}_{min}$-groups.
Thus it only leaves those of type $(2)$.
If $x\in L\setminus Z(G)$, then it was shown in the proof of Theorem~\ref{SolCAminGp} that $C_G(x)=C_L(x)$.
If $x \in G\setminus L$, then as $G/Z(G)$ is a Frobenius group $xZ(G)$ lies in some conjugate of $K/Z(G)$ \cite[Page 496]{Huppert}.
Thus assume $x\in K$.
As $K$ is abelian, then $K/Z(G)\leq C_G(x)/Z(G)\leq C_{G/Z(G)}(xZ(G))\leq K/Z(G)$.
That is $K=C_G(x)$.
It now follows that any two distinct minimal non-central element centralisers have intersection equal to $Z(G)$. 
\end{proof}
\end{prop}

Thus it only remains to study the non-solvable ${\rm CA}_{min}$-groups.

\begin{thm}\label{NonSolCAminGp}
Let $G$ be a non-solvable group.
Then $G$ is a ${\rm CA}_{min}$ group if and only if  $G/Z(G)\cong PSL_2(p^n)$ or $PGL_2(p^n)$ with $p^n>3$.
\begin{proof}
As $G/Z(G)$ admits a non-trivial normal partition, we will use the classification of Suzuki (Theorem~\ref{SuzNonSolPart}) to determine $G/Z(G)$.

\underline{\bf Case (1)}\newline
If $G/Z(G)$ is a Frobenius group, then as in the solvable case the kernel $L/Z(G)$ is nilpotent and the complement $K/Z(G)$ is abelian.
However, this implies $G/Z(G)$ and therefore $G$ is solvable.

\underline{\bf Case (2)}\newline
If $G/Z(G)$ is isomorphic to $Sz(2^n)$, then Schmidt[8] showed that the Sylow $2$-subgroups of $Sz(2^n)$ are subgroups of some components for any non-trivial partition. 
However $Sz(2^n)$ has non-abelian Sylow $2$-subgroups and therefore cannot be a ${\rm CA}_{min}$-group.

\underline{\bf Case (3)}\newline
Assume $G/Z(G)\cong PSL_2(p^n)$ or $PGL_2(p^n)$ with $p^n>3$.
In this case we want to show that any group arising in this way is a ${\rm CA}_{min}$-group.

It is well known that every element centraliser in $PGL_2(q)$ and $PSL_2(q)$ takes one of the following forms:
\begin{enumerate}
\item A cyclic group of order $q$, $q-1$ or $q+1$.
\item A dihedral group of order $2(q+1)$ or $2(q-1)$.
\end{enumerate}

Let $x\in G\setminus Z(G)$ such that $C_G(x)$ is a minimal centraliser and set $C$ to be the subgroup of $G$ such that $C_{G/Z(G)}(xZ(G))=C/Z(G)$.
If $C/Z(G)$ is cyclic, then $C$ is abelian and thus $C_G(x)=C$ .
If $C/Z(G)$ is dihedral, then as $xZ(G)\in Z(C/Z(G))$, it follows that $xZ(G)\in C'/Z(G)$ for $C'/Z(G)$ the cyclic subgroup of index $2$ in $C/Z(G)$.
Hence $C'\leq C_G(x)\leq C$.
If $C_G(x)=C$, then there exists a $y\in G$ such that $C_{G/Z(G)}(yZ)=C'/Z(G)$ and so $C_G(y)<C_G(x)$ contradicting minimality.
Thus every minimal centraliser in $G$ is abelian and its quotient is a centraliser in $G/Z(G)$.

Let $x,y\in G\setminus Z(G)$ such that $C_G(x)$ and $C_G(y)$ are distinct minimal centralisers in $G$.
Then $C_G(x)/Z(G)$ and $C_G(y)/Z(G)$ are centralisers in $G/Z(G)$.
It is enough to show that $(C_G(x)/Z(G))\cap (C_G(y)/Z(G))$ is trivial in $G/Z(G)$.
Let $kZ(G)\in (C_G(x)/Z(G)\cap C_G(y)/Z(G))$.
Then $C_G(x)/Z(G)$ and $C_G(y)/Z(G)$ are distinct abelian centralisers in $G/Z(G)$ which are subgroups of $C_{G/Z(G)}(kZ)$.
However, no centraliser in $PSL_2(q)$ or $PGL_2(q)$ contains two distinct abelian centralisers in $PSL_2(q)$ or $PGL_2(q)$ respectively.
Therefore $kZ=Z$ and hence the intersection is trivial.
In particular, we have shown that $C_G(x)\cap C_G(y)=Z(G)$ and $G$ is a ${\rm CA}_{min}$-group.
\end{proof}
\end{thm}

\begin{cor*}[Corollary~\ref{NonSolFIsCA}]
Any non-solvable F-group is a CA-group.
\begin{proof}
By Rebmann, any non-solvable F-group must be of the form $G/Z(G)\cong PSL_2(q)$ or $PGL_2(q)$ with some extra condition on $G'$.
However, by the previous result any group such that $G/Z(G)$ has this structure is ${\rm CA}_{min}$-group.
Hence $G$ is an F-group and ${\rm CA}_{min}$-group, and it follows it is a CA-group. 
\end{proof}
\end{cor*}

\subsection{Classifying ${\rm F}_{min}$-groups}

We now repeat a similar argument as in the case of ${\rm CA}_{min}$-groups. Most details are omitted, however we include details for cases (1) and (2) in the solvable case to highlight that the partition now consists of quotients of centres of centralisers. 

\begin{thm}
Let $G$ be a solvable ${\rm F}_{min}$-group.
Then $G$ is one of the following:
\begin{enumerate}
\item $G/Z(G)$ is a Frobenius group with kernel $L/Z(G)$ and complement $K/Z(G)$ such that both $K$ and $L$ are abelian.
\item $G/Z(G)$ is a Frobenius group with kernel $L/Z(G)$ and complement $K/Z(G)$ such that $K$ is abelian, $L$ is an ${\rm F}_{min}$-group, $Z(L)=Z(G)$ and $L/Z(L)$ is a $p$-group.
\item $G/Z(G)\cong {\rm Sym}(4)$ and if $V/Z(G)\cong V_4$, then $V$ is non-abelian. 
\item $G$ has an abelian normal subgroup of index $p$, $G$ is not abelian.
\item $G\cong A\times P$, where $A$ is abelian and $P$ is a non-abelian $p$-group for some prime $p$; therefore $P$ is an ${\rm F}_{min}$-group.
\end{enumerate}

\begin{proof}
As $G/Z(G)$ admits a non-trivial normal partition, we apply the classification of Baer (Theorem~\ref{BaerSolPart}) to determine $G/Z(G)$ and then $G$.
Note that cases (3) and (4) use the exact same argument as in Theorem~\ref{SolCAminGp} and so we shall not repeat them.

\underline{\bf Case (1)}\newline
Let $L/Z(G)$ denote the Frobenius kernel of $G/Z(G)$ and $K/Z(G)$ an element in the partition of $G/Z(G)$ which is self-normalising.
Then $K=Z(C_G(x))$ for some minimal centraliser $C_G(x)$ in $G$.
Using the same argument as in Theorem~\ref{SolCAminGp} shows $K/Z(G)$ is a Frobenius complement and abelian.
Furthermore, recall that for $x\in L\setminus Z(G)$, then $C_G(x)=C_L(x)$ and $C_L(x)$ is minimal centraliser in $L$ implies that $C_G(x)$ is minimal centraliser in $G$.

If $L$ has a unique minimal centraliser, then $L$ is abelian by Corollary~\ref{CAminCenAb} and thus of type $(1)$.
Hence assume $L$ has two distinct minimal centralisers $C_L(x)<L$ and $C_L(y)<L$ for $x,y\in L\setminus Z(G)$.
Then $C_L(x)=C_G(x)$ and $C_L(y)=C_G(x)$ are distinct minimal centralisers in $G$.
Moreover $Z(L)\leq Z(C_L(x))\cap Z(C_L(y))=Z(G)$ and so $Z(L)=Z(G)$.
In particular, we have shown that $L$ is an ${\rm F}_{min}$-group.
Repeating the argument in Theorem~\ref{SolCAminGp} also implies $G$ is of type $(2)$.

\underline{\bf Case (2)}\newline
In this case $G/Z(G)\cong {\rm Sym}(4)$. 
Let $V\leq G$ such that $V/Z(G)\cong V_4$.
If $V$ is abelian, then for all $x\in V\setminus Z(G)$ we have $V\leq C_G(x)$.
As $xZ(G)$ has order $2$ it follows that $C_{G/Z(G)}(xZ(G))=D_8$ or $V_4$ (when $xZ(G)$ is a double or single transposition respectively), we also observe that $V/Z(G)\leq C_G(x)/Z(G)\leq  C_{G/Z(G)}(xZ(G))$.

If there exists an $x\in V$ such that $C_G(x)/Z(G)\cong V_4$, then $C_G(x)=V$ and so is abelian and thus a minimal centraliser.
However $V_4$ is not contained in the partition of $G/Z(G)$ which yields a contradiction. 
Thus for all $x\in V\setminus Z(G)$ we have that $C_G(x)/Z(G)=C_{G/Z(G)}(xZ(G))\cong D_8$ and $xZ(G)\leq Z(D_8)$.
In particular, it follows that $V_4$ must be the normal Klein-four subgroup of ${\rm Sym}(4)$.

Inside $C_G(x)/Z(G)\cong D_8$ there exists a unique cyclic subgroup of order 4 and another copy of $V_4$ which is not normal in ${\rm Sym}(4)$.
Let $N,M\leq C_G(x)$ such that $N/Z(G)\cong C_4$ and $M/Z(G)$ equals the non-normal copy of $V_4$.
By Theorem~\ref{BaerSolPart}, $N/Z(G)$ lies in the partition and therefore $N$ is the centre of a minimal centraliser in $G$.
The subgroup $M$ contains $x$ and therefore $\langle x,Z(G)\rangle \leq Z(M)$. 
In particular, we see that $M/Z(M)$ must be cyclic and so $M$ is abelian.
However $x\in Z(N)\cap Z(M)$ and so $M$ cannot be a minimal centraliser in $G$.
Thus $C_G(y)>M$ for all $y\in M$.
As $M/Z(G)$ lies in a unique maximal subgroup isomorphic to $D_8$, it follows that for each $y\in M\setminus Z(G)$, then $C_G(y)/Z(G)$ equals the unique maximal subgroup containing $M/Z(G)$.
Thus $M/Z(G)\leq Z(D_8)$ which is a contradiction. 
\end{proof}

\end{thm}

Using the same arguments as in Proposition~\ref{ListAreCAmin} gives the analogous result for ${\rm F}_{min}$-groups.

\begin{prop}
Any solvable group occurring in Theorem~\ref{MainThm} is an ${\rm F}_{min}$-group.
\end{prop}

We are therefore left to study the non-solvable ${\rm F}_{min}$-groups.

\begin{thm}
Let $G$ be a non-solvable group.
Then $G$ is an ${\rm F}_{min}$-group if and only if  $G/Z\cong PSL_2(p^n)$ or $PGL_2(p^n)$ with $p^n>3$.
\begin{proof}
As $G/Z(G)$ admits a non-trivial normal partition, we will use the classification of Suzuki (Theorem~\ref{SuzNonSolPart}) to determine $G/Z(G)$.
Note that Cases (1) and (2) use exactly the same argument as in Theorem~\ref{NonSolCAminGp} and so shall not be repeated.

\underline{\bf Case (3)}\newline
Assume $G/Z(G)\cong PSL_2(p^n)$ or $PGL_2(p^n)$ with $p^n>3$.
We saw that any such group is a ${\rm CA}_{min}$-group, which implies it is an ${\rm F}_{min}$-group.
\end{proof}
\end{thm}

Moreover, we now obtain the analogous corollary of Rebmann for non-solvable ${\rm F}_{min}$-groups.

\begin{cor}
Any non-solvable ${\rm F}_{min}$-group is a ${\rm CA}_{min}$-group.
\end{cor}

Finally, by combing the two theorems in this section we obtain Theorem~\ref{MainThm}.

\section{A family of F-groups which are not CA-groups}
As we saw in the introduction, given the classification of ${\rm CA}_{min}$-groups, it is easy to see that there is a non-solvable ${\rm CA}_{min}$-group which is not a CA-group.
However, as commented in Rebmann, any non-solvable F-group is also a CA-group.
Thus we need to consider the solvable classification from Rebmann.
In particular, using \cite[Corollary 5.1]{FGroups} if $G$ is an F-group that is not a CA-group, $G$ must take one of the two forms:

\begin{enumerate}
\item $G\cong A\times P$ where $A$ is abelian and $P$ is a non-abelian $p$-group which is also an F-group.
\item $G/Z(G)$ is a Frobenius group with kernel $L/Z(G)$ and complement $K/Z(G)$, $K$ is abelian $Z(L)=Z(G)$, $L$ is an $F$-group and $L/Z(L)$ is a $p$-group. 
\end{enumerate}

Note that if we have an F-group which is not a CA-group of the second type, then the subgroup $L$ cannot be a CA-group.
In particular, if there exists an F-group which is not a CA-group, then there exists such a $p$-group.
Thus to find an F-group which is not a CA-group, the first place to consider is in the set of $p$-groups.
In particular we shall consider the class of extraspecial groups.

First we state the following lemma which will be of use to us.
\begin{lm}
Let $G$ be a finite group in which the derived subgroup has order $p$ for some prime $p$.
Then $G$ is an F-group.
\begin{proof}
This result follows from the observation that the conjugacy class $g^G$ is contained in the coset $gG'$.
Therefore $|g^G|$ equals $p$ or $1$.
Hence $|C_G(g)|=\frac{|G|}{p}$ or $|G|$ and so every non-central element centraliser is both maximal and minimal.
In particular $G$ is an F-group.
\end{proof}
\end{lm}

Let $G$ be an extraspecial group, usually denoted by one of the two groups $p^{1+2n}_{\pm}$ for some positive integer $n$.
Then we have $G'=\Phi(G)=Z(G)$ of order $p$.
By the previous lemma $G$ is an F-group.

Assume $n>1$, otherwise $G$ is of order $p^3$ and it is easy to see such groups are CA-groups.
Then $G$ is isomorphic to the central product of $H$ and $P$, where $H$ is an extraspecial group of order $p^{1\pm 2(n-1)}$ and $P$ is extraspecial of order $p^3$.
Take $x\in H\setminus Z(H)$.
Then $x\not\in Z(G)$, however $P\leq C_G(x)$.
Thus $G$ is not a CA-group.

\begin{prop*}[Proposition~\ref{FNotCA}]
Let $G$ be an extraspecial group of order $p^{2n+1}$ with $n>1$.
Then $G$ is an F-group which is not a CA-group.
\end{prop*}

Note that not all F-groups which are not CA-groups occur from extraspecial groups.
In particular, using GAP we can find 5 groups of order 64 which are F-groups but not CA-groups.

\section*{Acknowledgments}
The author gratefully acknowledges financial support by the ERC Advanced Grant $291512$.
In addition the author would like to thank Benjamin Sambale for reading and discussing a preliminary version of this paper.

\bibliographystyle{alpha}
\bibliography{bibfile}

\end{document}